\documentclass[11pt%
]{article}

\usepackage{amsmath}
\usepackage{amssymb}
\usepackage{amsthm}

\usepackage{graphics}
\usepackage{epigraph}
\usepackage{soul}

\usepackage{smart}
\SectioningParameter{presecnum}{section}=\noexpand\S;
\RestoreLaTeXeqno

\usepackage[all]{xy}

\usepackage{lscape}
\usepackage{rotating}

\usepackage{bbm}

\newcounter{notes}
\newenvironment{arab}{%
 \begin{list}{\arabic{notes})}{\usecounter{notes}%
  \settowidth{\labelwidth}{0.0em}
  \setlength{\labelsep}{0.0em}
  \setlength{\leftmargin}{0em}
  \setlength{\itemindent}{0.0em}
  \setlength{\itemsep}{1mm}
  \setlength{\topsep}{2mm}
  \setlength{\parsep}{3mm}
  \setlength{\partopsep}{0mm}
  }
 }
{\end{list}}

\newcommand {\supplus}{\mathop{{\supset}\llap{\raise
0.5pt\hbox{\normalfont\small+}\hskip 0.5pt}}}
\newcommand {\subplus}{\mathop{{\subset}\llap{\raise
0.5pt\hbox{\normalfont\small+}\hskip 0.5pt}}}
\newcommand {\divby}  {\lower 0.15ex \hbox{\,\vdots\,}}
\newcommand {\Cee}    {{\mathbb  C}}

\newcommand {\Zee}    {{\mathbb  Z}}

\newcommand {\fel}     {{\mathfrak{el}}}

\newcommand {\fg}     {{\mathfrak{g}}}    %

\newcommand {\fo}     {{\mathfrak{o}}}

\def \opname#1#2%
  {\expandafter\newcommand \csname #1\endcsname {{\mathop{\mathrm{#2}}\nolimits }}}

\newcommand{\rmname}[1]
  {\expandafter\newcommand \csname #1\endcsname {{\operatorname{#1}}}}

\newcommand{\rmnameii}[2]
  {\expandafter\newcommand \csname #1\endcsname {{\operatorname{#2}}}}

\rmname{act} \rmname{Ad} \rmname{Add} \rmname{ad} \rmname{Alt}
\rmname{alt} \rmname{Ann} \rmname{antidiag} \rmname{Ber}
\rmname{ber} \rmname{Bil} \rmname{Br} \rmname{card} \rmname{ch}
\rmname{Char} \rmname{cem} \rmname{cj} \rmname{Cliff}
\rmname{cntr} \rmname{codim} \rmname{Coind} \rmname{const}
\rmname{col} \rmname{cork} \rmname{cpr} \rmname{diag}
 \rmnameii{Div}{div} \rmname{Def} \rmname{Deg}
\rmname{Der} \rmname{Diff} \rmname{Dim} \rmname{End} \rmname{Even}
\rmname{Ext} \rmname{gr} \rmname{Hom} \rmname{HT}
\rmnameii{Ht}{ht} \rmname{hwt} \rmname{Id} \rmname{id}
\rmname{ind} \rmname{Ind} \rmname{Inf} \rmname{irr} \rmname{Le}
\rmname{Lie} \rmname{lwt} \rmname{mult} \rmname{Mat} \rmname{Mor}
\rmname{nm} \rmname{Ob} \rmname{Odd} \rmname{Osc} \rmname{per}
\rmname{Pic} \rmname{pr} \rmname{pro} \rmname{Prime} \rmname{Proj}
\rmname{prt} \rmname{pt} \rmname{Q} \rmname{qet} \rmname{qtr}
\rmname{rd} \rmname{rk} \rmname{row} \rmname{Res} \rmname{salt}
\rmname{Sch} \rmname{SBr} \rmname{sdim}\rmname{scalar}
\rmname{Ser} \rmname{sign} \rmname{Smbl} \rmname{spin}
\rmname{ssym} \rmname{str} %\rmname{st}
\rmname{sgn} \rmname{sq}
\rmname{symm} \rmname{supp} \rmname{Supp} \rmname{St}
\rmname{Spec} \rmname{Spm} \rmname{tr} \rmname{vpt} \rmname{Vect}
\rmname{weyl} \rmname{Weyl} \rmname{Witt}

\opname{vvol}  {{v\hspace{-0.1ex}o\hspace{-0.02ex}l\/}}
\opname{pnt}  {\text{\normalfont pt}} \opname{Span} {{Span}}
\opname{slim} {\overline{\lim}} \opname{Vol}
{{V\hspace{-0.55ex}o\hspace{-0.02ex}l\/}} \opname{Par}
{{P\hspace{-0.3ex}a\hspace{-0.05ex}r\/}}

\newcommand {\ev} {{\bar0}}
\newcommand {\od} {{\bar1}}

\newcommand {\tto} {\longrightarrow}

\newcommand {\bcdot}   {\mathbin{\hbox{\raise.4ex\hbox{\bf.}}}} % bold \cdot

\newcommand {\secno} {}

\makeatletter

\theoremstyle{plain}
\@@newtheorem*{Theorem}{\secno Theorem}
\@@newtheorem*{Lemma}{\secno Lemma}
\@@newtheorem*{Proposition}{\secno Proposition}
\@@newtheorem*{Corollary}{\secno Corollary}
\@@newtheorem*{Statement}{\secno Statement}
\@@newtheorem*{Problem}{\secno Problem}
\@@newtheorem*{Question}{\secno Question}
\@@newtheorem*{Conjecture}{\secno Conjecture}

\theoremstyle{definition}
\@@newtheorem*{Example}{\secno Example}
\@@newtheorem*{Examples}{\secno Examples}
\@@newtheorem*{Convention}{\secno Convention}
\@@newtheorem*{Comment}{\secno Comment}

\theoremstyle{remark}

\@@newtheorem*{Remark}{\secno Remark}
\@@newtheorem*{Remarks}{\secno Remarks}
\@@newtheorem*{Exercise}{\secno Exercise}
\@@newtheorem*{Solution}{\secno Solution}
\@@newtheorem*{Hint}{\secno Hint}

\newcommand{\?}{\nobreak\hskip.145em\nobreak\hskip\z@skip}

\renewcommand\appendix{\par
  \setcounter{section}{0}%
  \setcounter{subsection}{0}%
  \gdef\thesection{\appendixname\kern1ex\@Alph\c@section}}

\makeatother

\newcommand{\ssec}[2]{\subsection*{%
 \refstepcounter{subsection}%
 \boldmath\thesubsection.\kern1ex#2.}%
 \label{ss#1}%
}

\date{}

\begin{document}

\graphicspath{{./figs/}}

\title%[Cunha and Elduque superalgebras]
{Cartan matrices and presentations of the exceptional simple Elduque
Lie superalgebra}

\author{Sofiane Bouarroudj${}^1$, Pavel Grozman${}^2$, Dimitry Leites${}^3$
\thanks{DL is thankful to MPIMiS, Leipzig, for financial support and most
creative environment; to A.~Elduque for comments; to O.~Shirokova
for a \TeX pert help and to A.~Protopopov for a help with the
graphics, see \cite{Pro}.} }

\address{${}^1$Department of Mathematics, United Arab Emirates University, Al
Ain, PO. Box: 17551; Bouarroudj.sofiane@uaeu.ac.ae\\
${}^2$Equa Simulation AB, Stockholm, Sweden; pavel@rixtele.com\\
${}^3$MPIMiS, Inselstr. 22, DE-04103 Leipzig, Germany\\
on leave from Department of Mathematics, University of Stockholm,
Roslagsv. 101, Kr\"aft\-riket hus 6, SE-104 05 Stockholm,
Sweden; mleites@math.su.se, leites@mis.mpg.de}



\maketitle

\begin{abstract} Recently Alberto Elduque
listed all simple and graded modulo 2 finite dimensional Lie
algebras and superalgebras whose odd component is the spinor
representation of the orthogonal Lie algebra equal to the even
component, and discovered one exceptional such Lie superalgebra in
characteristic 5. For this Lie superalgebra all inequivalent
Cartan matrices (in other words, inequivalent systems of simple
roots) are listed together with defining relations between analogs
of its Chevalley generators.
\end{abstract}

%\date{Received May 1, 2006}

\section{Introduction}
In \cite{BGL}, for the ten exceptional simple finite dimensional
Lie superalgebras of Elduque and Cunha over an algebraically
closed field of characteristic 3, we have listed all inequivalent
Cartan matrices (in other words, inequivalent systems of simple
roots) and their inverse, if invertible, defining relations
between analogs of the Chevalley generators, and the coefficients
of linear dependence over $\Zee$ of the maximal roots with respect
to simple roots.

Here we provide with the same type of information for the only for
$p=5$ exceptional (new) simple finite dimensional Lie superalgebra
$\fel$ among the Lie superalgebras whose odd component is the
spinor representation of the orthogonal Lie algebra equal to the
even component. A.~Elduque described $\fel$ in components:
$\fel_\ev=\fo(11)$, $\fel_\od=\Gamma_5$ (in notations of
\cite{FH}, p. 376; i.e., the spinor representation), so $\sdim\
\fel =(55|32)$. To determine $\fel$, we have to explicitly define
the bracket: $[\cdot, \cdot]: S^2(\fel_\od)\tto \fel_\ev$, which
is not that easy; for details, see \cite{El2}. Since  the bracket
is not easy to describe in these terms, it was not clear if $\fel$
possesses a Cartan matrix (we found out that it does).

Not every simple finite dimensional Lie superalgebra $\fg$
possesses a Cartan matrix, but if it does, it is often more
convenient to express $\fg$ in terms of its Cartan matrices; for
their definition and description of the analogs of Dynkin graphs,
see \cite{GL1,BGL}. With the help of {\bf SuperLie} package
\cite{Gr} we list all inequivalent Cartan matrices (hence, systems
of simple roots) of $\fel$ and describe the defining relations
corresponding to these matrices.

For Lie superalgebras with Cartan matrix generated by Chevalley
generators, there are two types of defining relations: Serre-type
ones and non-Serre type ones (over $\Cee$, all these relations are
listed in \cite{GL1}). Sometimes some of the Serre-type relations
are redundant but this does not matter in practical calculations.
At the moment, the problem \so{how to encode the non-Serre type
relations in terms of Cartan matrix} is open for $p>0$. Some
relations (for any $p$) are so complicated that we conjecture that
there is no {\it general} encoding procedure. This is why our list
of relations is of practical interest.

For the background on Cartan matrices and \lq\lq odd reflections"
that connect inequivalent matrices, see \cite{BGL}.

\section{The Elduque superalgebra: Systems of simple roots}

Having found out one Cartan matrix of $\fel$, we list them all. We
denote by $\fel^{n)}$ the realization of $\fel$ by means of the
$n$th Cartan matrix.

The table 8) shows the result of odd reflections (the number $n$
of the row is the number of the matrix $n)$ in the list below, the
number of the column is the the number of the root in which
reflection is made; the cells contain the results of reflections
(the number of the matrix obtained) or a \lq\lq --" if the
reflection is not appropriate because $A_{ii}\neq 0$. The nodes
are numbered by small boxed numbers; instead of joining nodes with
four segments in the cases where $A_{ij}=A_{ji}=1\equiv -4\mod 5$
we use one dotted segment; curly lines with arrows depict odd
reflections.

\begin{figure}[ht]%\centering
\parbox{.7\linewidth}{\includegraphics{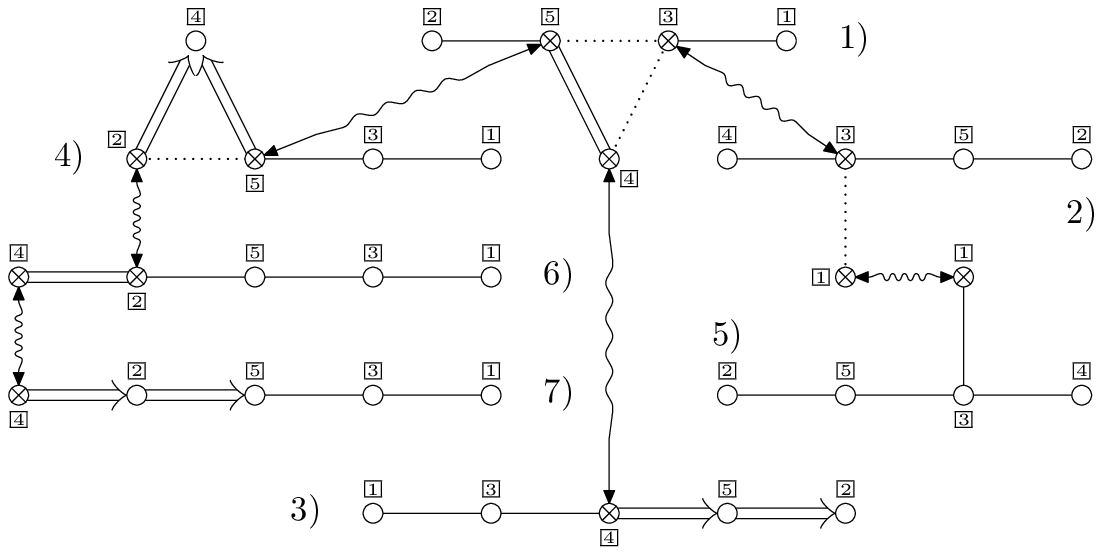}}
\end{figure}

{\tiny %footnotesize

$$
\begin{matrix}
1) \begin{pmatrix}
2&0&-1&0&0 \\
    0&2&0&0&-1 \\
    -1&0&0&1&1 \\
    0&0&1&0&-2 \\
    0&-1&1&-2&0
\end{pmatrix}\quad
2) \begin{pmatrix}
    0&0&1&0&0 \\
    0&2&0&0&-1 \\
    1&0&0&-1&-1 \\
    0&0&-1&2&0 \\
    0&-1&-1&0&2
    \end{pmatrix}
\quad
    \boxed{3)} \begin{pmatrix}
    2&0&-1&0&0 \\
    0&2&0&0&-1 \\
    -1&0&2&-1&0 \\
    0&0&-1&0&2 \\
    0&-2&0&-1&2
    \end{pmatrix}\end{matrix}
$$
$$
\begin{matrix}
4) \begin{pmatrix}
    2&0&-1&0&0 \\
    0&0&0&2&1 \\
    -1&0&2&0&-1 \\
    0&-1&0&2&-1 \\
    0&1&-1&2&0
    \end{pmatrix}
\quad \boxed{5)} \begin{pmatrix}
    0&0&-1&0&0 \\
    0&2&0&0&-1 \\
    -1&0&2&-1&-1 \\
    0&0&-1&2&0 \\
    0&-1&-1&0&2
    \end{pmatrix}\quad
6) \begin{pmatrix}
    2&0&-1&0&0 \\
    0&0&0&-2&-1 \\
    -1&0&2&0&-1 \\
    0&-2&0&0&0 \\
    0&-1&-1&0&2
    \end{pmatrix}
\end{matrix}
$$
$$
\begin{matrix}
    \boxed{7)} \begin{pmatrix}
    2&0&-1&0&0 \\
    0&2&0&-1&-2 \\
    -1&0&2&0&-1 \\
    0&2&0&0&0 \\
    0&-1&-1&0&2
    \end{pmatrix}\quad
    8) \begin{pmatrix} %{\footnotesize
    -&-&2&3&4 \\
    5&-&1&-&- \\
    -&-&-&1&- \\
    -&6&-&-&1 \\
    2&-&-&-&- \\
    -&4&-&7&- \\
    -&-&-&6&-
  \end{pmatrix}
\end{matrix}
$$
}

\medskip

\normalsize

\section{The Elduque superalgebra: Defining relations}
{\bf To save space, we omit indicating the Serre relations in what
follows; their fulfilment is assumed}.

{\footnotesize
\begin{arab}

\item[]

\item[]
$
\arraycolsep=0pt
\renewcommand{\arraystretch}{1.2}
\begin{array}{ll}
1)\  &
  {{\left[x_{4},\,\left[x_{3},\,x_{5}\right]\right]}-
  {\left[x_{5},\,\left[x_{3},\,x_{4}\right]\right]}=0};\quad
  {{\left[\left[x_{1},\,x_{3}\right],\,\left[x_{3},\,x_{4}\right]\right]}=
  0};\quad
  {{\left[\left[x_{1},\,x_{3}\right],\,\left[x_{3},\,x_{5}\right]\right]}=
  0};\\
&
{{\left[\left[x_{2},\,x_{5}\right],\,\left[x_{3},\,x_{5}\right]\right]}=
  0};\quad
  {{\left[\left[x_{4},\,x_{5}\right],\,\left[\left[x_{2},\,x_{5}\right],\,
  \left[x_{4},\,x_{5}\right]\right]\right]}
     = 0}\\
2) \ &
  {{\left[\left[x_{1},\,x_{3}\right],\,\left[x_{3},\,x_{4}\right]\right]}=
  0};\quad
  {{\left[\left[x_{1},\,x_{3}\right],\,\left[x_{3},\,x_{5}\right]\right]}=
  0};\\
  &
  {{\left[\left[\left[x_{3},\,x_{4}\right],\,\left[x_{3},\,x_{5}\right]\right],\,
      \left[\left[x_{3},\,\left[x_{2},\,x_{5}\right]\right],\,
       \left[\left[x_{3},\,x_{4}\right],\,\left[x_{3},\,x_{5}\right]\right]\right]\right]}=
       0}\\
3) \ & {{\left[\left[x_{3},\,x_{4}\right],\,\left[\left[x_{2},\,x_{5}\right],\,
\left[x_{4},\,x_{5}\right]\right]\right]}
     -{4\, \left[\left[x_{4},\,\left[x_{2},\,x_{5}\right]\right],\,
\left[x_{5},\,\left[x_{3},\,x_{4}\right]\right]
        \right]}=0};\\
&        {{\left[\left[x_{4},\,\left[x_{1},\,x_{3}\right]\right],\,
      \left[\left[x_{3},\,x_{4}\right],\,\left[x_{4},\,x_{5}\right]\right]\right]}=0}\\
4)\ &
{{\left[x_{4},\,\left[x_{2},\,x_{5}\right]\right]}-
{3\, \left[x_{5},\,\left[x_{2},\,x_{4}\right]\right]}=0}
  ;\quad
  {{\left[\left[x_{2},\,x_{5}\right],\,\left[x_{3},\,x_{5}\right]\right]}=
  0};\\
  &
  {{\left[\left[x_{5},\,\left[x_{1},\,x_{3}\right]\right],\,
      \left[\left[x_{3},\,x_{5}\right],\,\left[x_{4},\,x_{5}\right]\right]\right]}=0}\\
5)\ &
{{\left[\left[\left[x_{4},\,\left[x_{1},\,x_{3}\right]\right],\,
  \left[x_{5},\,\left[x_{1},\,x_{3}\right]\right]\right],\,
      \left[\left[\left[x_{1},\,x_{3}\right],\,\left[x_{2},\,x_{5}\right]\right],\,
       \left[\left[x_{4},\,\left[x_{1},\,x_{3}\right]\right],\,
       \left[x_{5},\,\left[x_{1},\,x_{3}\right]\right]\right]\right]
       \right]}= 0}\\
6)\ &
  {{\left[\left[x_{2},\,x_{4}\right],\,\left[\left[x_{2},\,x_{4}\right],\,
\left[x_{2},\,x_{5}\right]\right]\right]}
     = 0};\\
     & {{\left[\left[\left[x_{1},\,x_{3}\right],\,\left[x_{2},\,x_{5}\right]\right],\,
      \left[\left[x_{3},\,\left[x_{2},\,x_{5}\right]\right],\,
       \left[\left[x_{2},\,x_{4}\right],\,\left[x_{2},\,
x_{5}\right]\right]\right]\right]}=0}\\
7)\ &
  {{\left[x_{2},\,\left[x_{2},\,\left[x_{2},\,x_{5}\right]\right]\right]}=
  0};\\
  &
  {\left[\left[\left[x_{2},\,x_{4}\right],\,\left[x_{5},\,\left[x_{1},\,
x_{3}\right]\right]\right],\,
      \left[\left[x_{5},\,\left[x_{2},\,x_{4}\right]\right],\,
       \left[x_{3},\,\left[x_{2},\,\left[x_{2},\,x_{5}\right]\right]\right]\right]\right]}\\
       &
    -{2\, \left[\left[\left[x_{2},\,x_{4}\right],\,
\left[\left[x_{1},\,x_{3}\right],\,\left[x_{2},\,x_{5}\right]\right]\right]
        ,\,\left[\left[x_{3},\,\left[x_{2},\,x_{5}\right]\right],\,
\left[x_{5},\,\left[x_{2},\,x_{4}\right]\right]\right]
        \right]}=0
      \end{array}
$ \end{arab}
}

 \ssec{}{The maximal roots} The maximal roots are of hight 41 and
their weights are equal to $(1,0,0,0,0)$, except root 5) whose
weight is equal to $(4,0,0,0,0)$. The coefficients of their
decomposition with respect to simple roots (over $\Zee$) are as
follows:
\[
\arraycolsep=0pt
\renewcommand{\arraystretch}{1.1}
\begin{array}{llllllll}
 1)\ & (2, 2, 3, 3, 4),\;&
 2)\ &(2, 2, 6, 3, 4),\; &
 3)\ & (2, 2,  3, 4, 4),\;&
 4)\ & (2, 2, 3, 3, 4),\\
 5) \ & (5, 2, 6, 3, 4),\;&
 6) \ &  (2, 5, 3, 3, 4),\;&
 7) \ &  (2, 5, 3, 2, 4)&&
    \end{array}
    \]
%    \end{landscape}

%\newpage

\section{The inverse matrices of the Cartan matrices}
(Their numbering matches that of the Cartan matrices.)

{\tiny %footnotesize
$$
1)\begin{pmatrix} 2& 2& 3& 3&4\\
2& 4& 4& 0&2\\
3&4 &1 &1 &3\\
3&0 &1 &3 &0 \\
4& 2& 3 &0 &4
\end{pmatrix}\quad 2) \begin{pmatrix}
 2 & 2& 1& 3& 4\\
  2& 4& 0& 0&2\\
1& 0& 0& 0& 0\\
 3& 0& 0& 3&0\\
4& 2& 0& 0&4
\end{pmatrix}\quad
3)\begin{pmatrix} 2& 2& 3& 4& 2\\
2& 4& 4& 1& 1\\
3& 4& 1& 3& 4\\
 4& 1& 3& 2&1\\
  4& 2& 3& 2&2
\end{pmatrix}
\quad 4)\begin{pmatrix}
2 & 2& 3& 4& 4\\
 2& 4& 4& 0& 1\\
  3& 4& 1& 3& 3\\
   3& 0& 1& 4&4\\
4&1&3&2&2
\end{pmatrix}
$$
$$ 5)\begin{pmatrix}
 0 & 3& 4& 2&1\\
 3&4&0&0&2\\
4& 0& 0& 0&0\\
 2&0&0&3&0\\
 1&2&0&0&4
\end{pmatrix}\quad 6) \begin{pmatrix}
2 & 0& 3& 3& 4\\
 0 &0& 0& 2& 0\\
  3& 0& 1& 1& 3\\
3& 2& 1& 3&4\\
4& 0& 3& 4& 2\\
\end{pmatrix}\quad
7) \begin{pmatrix}
2& 0& 3& 2&4\\
 0& 0& 0& 3&0\\
3& 0& 1& 4& 3\\
 2& 4& 4& 4&1\\
 4& 0& 3& 1&2\\
\end{pmatrix}
$$
}

\end{document}